\documentclass[10pt,leqno]{article}

\usepackage[francais,english]{babel}
\usepackage{amsmath,amssymb,amscd}
\usepackage{a4wide}
\usepackage[latin1]{inputenc}
\usepackage[all]{xy} 
\usepackage{hyperref}

\title{On the fine expansion of the unipotent contribution of the Guo-Jacquet trace formula}
\author{Pierre-Henri Chaudouard}
\date{}

\makeatletter
\@addtoreset{equation}{subsection}         
\makeatother

\newenvironment{paragr}[1][]{\refstepcounter{subsubsection} \noindent \textbf{\thesubsubsection . \ #1}}{\medskip}

\newenvironment{theoreme}{ \medskip\refstepcounter{theo}  \noindent\textbf{Theorem \thetheo}. ---\em}{\em \medskip}
\newenvironment{proposition}{\medskip\refstepcounter{theo}   \noindent\textbf{Proposition \thetheo}. ---\em}{\em\medskip}
\newenvironment{corollaire}{\medskip\refstepcounter{theo}  \noindent\textbf{Corollary \thetheo}. ---\em}{\em\medskip}

\newenvironment{lemme}{\medskip\refstepcounter{theo}   \noindent\textbf{Lemma \thetheo}. ---\em}{\em\medskip}

\newenvironment{preuve}[1][]{\noindent \textbf{Proof.} #1 --- }{\hfill
  \ensuremath{\square} \medskip}

\newenvironment{remarque}{\medskip\refstepcounter{theo}  \noindent\textbf{Remark \thetheo}. ---}{\medskip}

\DeclareMathOperator{\vol}{vol}

\DeclareMathOperator{\Ad}{Ad}

\DeclareMathOperator{\End}{End}

\DeclareMathOperator{\Norm}{Norm}

\DeclareMathOperator{\Aut}{Aut}
\DeclareMathOperator{\Gal}{Gal}
\DeclareMathOperator{\Hom}{Hom}

\DeclareMathOperator{\Id}{Id}

\DeclareMathOperator{\Res}{Res}

\DeclareMathOperator{\trace}{trace}

\newcommand{\RR}{\mathbb{R}}
\newcommand{\AAA}{\mathbb{A}}
\newcommand{\CC}{\mathbb{C}}

\newcommand{\QQ}{\mathbb{Q}}

\newcommand{\oc}{\mathcal{O}}

\newcommand{\Sc}{\mathcal{S}}

\newcommand{\uc}{\mathcal{U}}
\newcommand{\vc}{\mathcal{V}}

\newcommand{\pc}{\mathcal{P}}

\newcommand{\nc}{\mathcal{N}}

\newcommand{\bc}{\mathcal{B}}

\newcommand{\mgo}{\mathfrak{m}}
\newcommand{\ngo}{\mathfrak{n}}

\newcommand{\cgo}{\mathfrak{c}}
\newcommand{\pgo}{\mathfrak{p}}

\newcommand{\sgo}{\mathfrak{s}}

\newcommand{\hgo}{\mathfrak{h}}

\newcommand{\ugo}{\mathfrak{u}}

\newcommand{\al}{\alpha}

\newcommand{\om}{\omega}
\newcommand{\Om}{\Omega}

\newcommand{\back}{\backslash}

\newcommand{\Cc}{C_c^\infty}

\newcommand{\bg}{\langle}
\newcommand{\bd}{\rangle}

\newcommand{\eps}{\varepsilon}

\renewcommand{\leq}{\leqslant}
\renewcommand{\geq}{\geqslant}

\begin{document}
\maketitle
\selectlanguage{english}

\begin{abstract}
  For a useful  class of functions (containing functions whose one finite component is essentially a matrix coefficient of  a supercuspidal representation), we establish three results about the unipotent contribution of the Guo-Jacquet relative trace formula for the pair $(GL_n(D),GL_n(E))$. First we get a fine expansion in terms of global nilpotent integrals. Second we express these  nilpotent integrals in terms of zeta integrals. Finally we prove that they satisfy certain homogeneity properties. The proof is based on a new kind of truncation introduced in a previous article.
\end{abstract}
\tableofcontents

\section{Introduction}

\subsection{The statement}

\begin{paragr}
  Let $E/F$ be a quadratic extension of number fields. Let $\eps\in E^\times$ such that $\trace_{E/F}(\eps)=0$. Let $n\geq 1$ be  an integer. Let  $G=GL(2n,F)$ and  $H=GL(n,E)$ viewed as algebraic group over $F$. The $F$-basis $(1,\eps)$ identifies $E^n$ with $F^{2n}$ and gives an embedding $H\subset G$.  We identify $\eps$ with a scalar matrix in $GL(n,E)$.  Then $H$ is the centralizer of $\eps$ in $G$. Moreover,  the quotient $G/H$ is a symmetric space which is easily identified with  the subvariety $S\subset G$ of $g\in G$ such that $\eps^{-1}g\eps=g^{-1}$. The group $H$ acts on $S$ by conjugation.
\end{paragr}

\begin{paragr} Let $\AAA$ be the ring of adèles of $F$ and $|\cdot|$ be the product of local normalized absolute values.  We identify $\eps$ with a central element in $GL(n,E)$.   For any cuspidal automorphic form $\varphi$ on $G(\AAA)$, one can define its $H$-period as
  $$\pc_H(\varphi)=\int_{H(F)\back H(\AAA)^1}\varphi(h)\,dh$$
  where $H(\AAA)^1\subset H(\AAA)$ is the kernel of the morphism $h\mapsto |\det(h)|$ and $dh$ is an invariant measure. A cuspidal automorphic representation is said to be $H$-distinguished if $\pc_H$ induces a non-zero linear form on its underlying space. A fundamental question is to understand the interplay between the $H$-distinction, the Jacquet-Langlands functoriality  and distinction by some other related subgroups: a beautiful answer is given by  the so-called Guo-Jacquet conjecture (cf.   conjecture in \cite{GuoCanad}  extrapolating on results of Waldspurger and Jacquet).
\end{paragr}

\begin{paragr}
  A promising tool to sudy this problem is the so-called Guo-Jacquet trace formula (a specific example of a relative trace formula) based on the seminal work of Jacquet (cf. \cite{Jac}). A simple form of it has already been used successfully (see \cite{FMW}). It consists in expressing geometrically and spectrally the integral
  \begin{align}\label{eq:FTR}
    \int_{H(F)\back H(\AAA)^1} \sum_{\gamma \in S(F)}\Phi(h^{-1}\gamma h)\, dh
  \end{align}
  where $\Phi\in \Cc(S(\AAA))$.
  In general, \eqref{eq:FTR} is not convergent. To remedy this problem, the simple trace formula introduces severe  restrictions at two places on the function $\Phi$ so that on the spectral side only the cuspidal spectrum appears and on the geometric side only regular  semisimple $\gamma$'s contribute to the rational sum in  \eqref{eq:FTR}. However, for the purpose of convergence, it suffices to impose a certain local condition at one place on $\Phi$.   When it is satisfied, we shall say that $\Phi$ is very cuspidal. In the paper, this is the local condition \eqref{eq:vanish} below. Here it suffices to say that the condition is satisfied if the local component of $\Phi$ at a finite place $v$ is obtained by integration of a matrix coefficient of a supercuspidal representation of $G(F_v)/Z(F_v)$ (see remark \eqref{rq:mat-coef}, $Z$ is the center of $G$).  

For a very cuspidal $\Phi$, as a part of \eqref{eq:FTR} we have the (convergent) unipotent contribution defined by
   \begin{align}\label{eq:FTRu}
    \int_{H(F)\back H(\AAA)^1} \sum_{X\in \nc(F)}\Phi(h^{-1}\exp(X)h)\, dh
  \end{align}
  where $\nc(F)$ is the cone of nilpotent matrices in the tangent space  $\sgo$ of $S$ at identity (given by matrices $X$ of size $2n$ such that $\eps X+X\eps=0$) and $\exp$ is the usual exponential. We emphasize that besides regular semi-simple terms the unipotent contribution
is certainly the most important contribution of the geometric side of the Guo-Jacquet trace formula to understand (for general terms there should be some kind of Jordan decomposition).

 \end{paragr}

\begin{paragr} The goal of the paper is to prove  the following properties  of  the unipotent contribution for a very cuspidal $\Phi$.

\begin{theoreme}\label{thm:intro}(for a more precise statement see theorem \ref{thm:homog-gp} below)
Let   $\Phi\in \Cc(G(\AAA))$ be  a very  cuspidal function. Let  $\oc\in \nc(F)/H(F)$ be a nilpotent orbit and let $v$ be a place of $F$. 
\begin{enumerate}
\item For any $t\in F_v^\times$,  the integral 
$$J_\oc^t(\Phi)=\int_{[H]^1} \sum_{X\in \oc} \Phi(h^{-1}\exp(t^{-1}X)h )\,dh$$
is absolutely convergent. In particular we get $J_\oc^1(\Phi)$ for $t=1$.
\item We have
$$  \int_{H(F)\back H(\AAA)^1}  \sum_{X\in \nc(F)}\Phi(h^{-1}\exp(X)h)\, dh=\sum_{\oc\in \nc(F)/H(F)} J_\oc^1(\Phi).$$
\item There exists a  bound $\eta >0$ such that for any $t$ in $F_v^\times$ such that $|t|_v<\eta$, we have
$$ J_\oc^t(\Phi)=\lim_{s\to 0^+} s \theta_\oc(s)  Z_\oc(\Phi,s)$$
where
\begin{itemize}
\item $\theta_\oc(s)$ is some holomorphic function which is defined for $s\in \CC$ such that its real part satisfies $\Re(s)>0$ and which does not depend on $\Phi$ 
\item the function $Z_\oc(\Phi,s)$ is a zeta function attached to $\Phi$ of the variable $s\in \CC$, holomorphic for  $\Re(s)>0$.
\end{itemize}
\item Let $\eta$ be the above bound. For any  $t_0,t $ in $F_v^\times$ such that $|t_0|_v<\eta$ and  $|t|_v\leq 1$  we have
     \begin{align*}
       J_\oc^{tt_0}(\Phi)= |t|_v^{\dim(\oc)/2}J_{\oc}^{t_0}(\Phi).
     \end{align*}
   \end{enumerate}
 \end{theoreme}

 \begin{remarque} Assertions 1 and 2 provide a ``fine expansion'' of the unipotent contribution that is an expansion according conjugacy classes. Assertion 3 gives a way to compute each nilpotent contribution; moreover the zeta function  $Z_\oc(\Phi,s)$ admits a Eulerian product so it is possible to give an expression of each nilpotent contribution in terms of local objects.
 \end{remarque}

\begin{remarque}
  In the paper, the theorem is in fact stated and proved in the broader situation where $G=GL_n(D)$ where $D$ is a quaternion algebra containing $E$. It should also hold in the situation where $G$ is the multiplicative group of a $F$-simple central algebra containing $E$ and $H$ is the centralizer of $E$ in $G$. The methods of the article should also apply in this context but we have not written the details.
\end{remarque}

We have the following corollary.

\begin{corollaire}\label{cor:intro}(see corollary  \ref{cor:homog-gp}) Let   $\Phi\in \Cc(G(\AAA))$ be  a very  cuspidal function. Let  $v$ be a place of $F$. 
When $t\in F_v^\times $ goes to $0$, the expression 
$$\int_{[H]^1} \sum_{X\in \nc(F)} \Phi(h^{-1}\exp(t^{-1}X)h )\, dh
$$
is equivalent to 
$$\vol([H]^1) \Phi(1).$$
\end{corollaire}

\begin{remarque}
  One of the motivations for this result is that it plays a role  in Xue's approach of local distinction problems. For example, corollary \ref{cor:intro} is used to prove the existence of $H$-distinguished cuspidal automorphic representations of $G$ with local supercuspidal components at some places (see corollary 6.2 of \cite{Xue}). An other potential application of theorem \ref{thm:intro} and its corollary is the obtention of some kind of (relative) Weyl's law.
\end{remarque}
\end{paragr}

\subsection{The methods}

\begin{paragr}
The Guo-Jacquet trace formula in our context should share many properties with the usual Arthur-Selberg trace formula. Recall that the unipotent contribution of the Arthur-Selberg trace formula has a fine expansion in terms of local unipotent integrals (although some global constants remain unknown; for some progress on these questions see \cite{scuft} and \cite{cuft}). As a consequence it does satisfy general homogeneity properties (see \cite{A-unip}). We could have tried to prove analogs of several results of Arthur in our context (mainly the results of \cite{A-TF1}, \cite{A-woi} and \cite{A-unip}). In this paper  we choose to offer a rather simple proof of   theorem \ref{thm:intro} based on ideas developped in \cite{suvta}. The payoff is that we get a (in principle) computable expression (rather than unknown global coefficients).
\end{paragr}

\begin{paragr}
  First one reduces the problem to the case of the infinitesimal situation where $H$ acts on the tangent space $\sgo$. This reduction is the content of the final section, section \ref{sec:unip}. We introduce in  §\ref{S:weak-cusp} the set of weakly cuspidal functions $f \in \Cc(\sgo(\AAA))$ that satisfy some vanishing conditions. For such functions, we have
  \begin{align}\label{eq:expan}
    \int_{[H]^1} \sum_{X\in \nc(F)} f(\Ad(h^{-1})X) \,dh =\sum_{\oc\in \nc(F)/H(F)}  \int_{[H]^1} \sum_{X\in \oc} f(\Ad(h^{-1})X) \,dh  ;
  \end{align}
 on the right-hand side we sum   terms  that satisfy the  non-trivial convergence result (see theorem \ref{cor:cv-inf}) 
  \begin{align*}
 \int_{[H]^1} |\sum_{X\in \oc} f(\Ad(h^{-1})X) |\,dh <\infty.
  \end{align*}
 However these terms are difficult to compute since one cannot permute the sum over $X$ and the integral. This is where enters the truncation of \cite{suvta}. It enables to recover each term as a limit
  \begin{align*}
  \lim_{s\to 0}     s J_\oc(s,f)
\end{align*}
where $J_\oc(s,f)$ is a holomorphic function   for $s\in \CC$ of real part $\Re(s)>0$ (see theorem \ref{thm:limit}). The point is that $J_\oc(s,f)$ can be expressed in terms of a zeta integral under a mild assumption on the support of $f$ (see theorem \ref{thm:comput}). For this zeta integral, the homogeneity property is easy to check (see lemma \ref{lem:homog}).
\end{paragr}

\subsection{Acknowledgement}

We thank Hang Xue for raising the problem. We thank also Huajie Li for many useful discussions.

\section{Infinitesimal situation}\label{sec:inf}

\subsection{Notations}\label{ssec:algprel}

\begin{paragr}
Let $F$ be a number field and $\tau\in F$ such that $E=F[\sqrt{\tau}]$ is a quadratic extension. Let $\sigma$ be the generator of the Galois group $\Gal(E/F)$. For each place $v$ of $F$, let $|\cdot|_v$ be the normalized absolute value.

Let $\AAA$ be the ring of adèles of $F$. Let $\AAA^\times$ be the multiplicative group of $\AAA$. For any $x=(x_v)_v\in \AAA^\times$, we shall denote by $|x|$ the product $\prod_v |x|_v$ over all places of $F$.

For any algebraic group $G$ over $F$, We shall write $[G]$ for the quotient $G(F)\back G(\AAA)$. 
\end{paragr}

\begin{paragr}
  Let $V_F$ be of a vector space of dimension $n$ over $F$. Let  $V=V_F\otimes_F E$. By abuse of notation, we denote by $\sigma$ the $F$-automorphism of $V$ given by  $\Id_{V_F}\otimes_F\sigma$. Let $\sgo\subset \End_F(V)$ be the $F$-subspace of $\sigma$-linear endomorphisms of $V$ namely the space of maps $X$ of $V$ into itself such that $X\circ\sigma$ is $E$-linear.
\end{paragr}

\begin{paragr}
  Let $H$ be the algebraic group of automorphisms of the $E$-vector space $V$. By restriction of scalars, we will view $H$ as an $F$-group.  The group $H$ is provided with a Galois $F$-automorphism still denoted by $\sigma$. 
The group $H$ acts by conjugation on $\End_F(V)$ : we denote by $\Ad$ the restriction of this action on $\sgo$. This action is defined over $F$.  By differentiation, we get an action $\hgo\times \sgo \to \sgo$ denoted by $(X,Y)\mapsto [X,Y]$. For any $X\in \sgo$, let $H_X$ be the centralizer of $X$ in $H$.
\end{paragr}

\begin{paragr}\label{S:det}
  By abuse of notation, we will simply  denote by $\det$ the $F$-morphism $H\to \mathbb{G}_{m,F}$ given by $\Norm_{E/F}\circ \det$.
Let $H(\AAA)^1$ be the kernel of $h\mapsto |\det(h)|$. We shall denote by $[H]^1$ the quotient $H(F)\back H(\AAA)^1$. 
\end{paragr}

\begin{paragr}[Parabolic decomposition.] --- \label{S:parab} We fix an ordered $F$-basis $(e_1,\ldots,e_n)$ of $V_F$. By a (standard) parabolic subgroup $P$, we mean a subgroup which stabilizes an incomplete  standard flag  of $E$-subspaces 
$$V_0=(0)\subsetneq V_1\subsetneq V_1\oplus V_2\subsetneq \ldots \subsetneq    V_1\oplus\cdots \oplus V_r
$$
where for $1\leq i\leq r$ the $E$-vector space $V_i$ is  generated by vectors $e_{d_i+1},e_{d_i+2},\ldots, e_{d_i+\dim(V_i)}$ where $d_i=\sum_{j=1}^{i-1}\dim(V_j)$.

Then we have a standard Levi decomposition $P=MN$ where $N$ is the unipotent radical of $P$ and $M$ is the common stabilizer of the subspaces $V_i$. We define $\sgo_M$, resp. $\sgo_N$, to be the subspace of $X\in \sgo$ such that $XV_i\subset V_i$, resp. $X(V_1\oplus \cdots \oplus V_i)\subset V_1\oplus \cdots \oplus V_{i-1}$ for all $1\leq i\leq r$.
  Let $\sgo_P=\sgo_M\oplus\sgo_N$.
  The groups $P, M$ and $N$ act respectively on $\sgo_P, \sgo_M$ and $\sgo_N$.

The stabilizer of the complete standard flag is denoted by $B$. We have a standard Levi decomposition $B=TN_B$.
\end{paragr}

\begin{paragr}[Maximal compact subgroup.] --- \label{S:maxcpt}Thanks to the basis of §\ref{S:parab},  we identify $H(\AAA)$ with $GL(n,\AAA_E)$. Let 
 $K\subset H(\AAA)$ be the maximal  compact subgroup corresponding to the standard maximal compact subgroup of $GL(n,\AAA_E)$. For parabolic subgroups $P=MN$, we have the Iwasawa decomposition $H(\AAA)=M(\AAA)N(\AAA)K$.  
\end{paragr}

\begin{paragr}[Haar measures.] ---  \label{S:Haar} We fix a Haar measure on $\AAA_E^\times$ the multiplicative group of the adèles of $E$. For any $n\geq 1$, we take on $(\AAA_E^\times)^n$ the product of the Haar measure. On any unipotent group $N$ over $F$, we take  the Haar measure on $N(\AAA)$ such that the quotient measure on $[N]$ gives the total volume $1$. The Haar measure on the standard compact maximal of $GL(n,\AAA_E)$ is such that the total volume is $1$. Finally we require that the Haar measure on  $GL(n,\AAA_E)$ is compatible with the Iwasawa decomposition, which normalizes the measure if we view $(\AAA_E^\times)^n$ as a minimal Levi subgroup. In this way, we get a normalization of the  Haar measure on  $H(\AAA)$. The Haar measure on $H(\AAA)^1$ is such that the quotient measure on  $H(\AAA)/H(\AAA)^1\simeq \RR_+^\times$ is the usual Haar measure. Then  we get quotient measures on $[H]$ and $[H]^1$.
\end{paragr}

\begin{paragr}[Categorical quotient.] ---  Let $C: \sgo \to \cgo=\sgo//H$ be the categorical quotient. The quotient $\cgo$ can be identified with the standard $n$-dimensional affine space over $F$ in such a way that $c$ is the map that associates to any $X\in \sgo$ the coefficients of the characteristic polynomial of $X^2\in \End_E(V)$ (the coefficients are in fact defined over $F$).  For any  $c\in \cgo$, let $\sgo_c$ be the fiber of $C$ above $c$. In particular, if $c$ corresponds to the polynomial $t^{n}$, we denote by $\nc$ the  fiber $\sgo_c$.  Then an element $X\in\sgo$ belongs to $\nc$ if and only if the $E$-endomorphism $X^2$ of $V$ is nilpotent in the usual sense. By abuse, in the following, $\nc$ will be  called the nilpotent cone and the elements of $\nc$ will be called the nilpotent elements of $\sgo$.
\end{paragr}

\subsection{Induction of nilpotent orbits}

\begin{paragr}[Classification of nilpotent orbits.] --- It is given by the following lemma.

  \begin{lemme}\label{lem:classif-orbit}
    For any $X\in \nc$, there exists an $E$-basis of $V$ such that the matrix of $X$ is in the Jordan normal form.

    There are finitely many orbits of $H$ on $\nc$ and they are classified by their  Jordan normal form.
      \end{lemme}

      \begin{preuve}
        One can find a proof in \cite{GuoPJ} lemma 2.3. Alternatively, one can observe that $V$ has a $E[t]$-module structure of $t^{2n}$-torsion given by $(P,v)\mapsto P(X)v$  which boils down to the classification of torsion modules over a principal ideal domain.
      \end{preuve}
\end{paragr}

\begin{paragr}[Induction à la Lusztig-Spaltenstein] ---  We state and prove in our context some results about induction of nilpotent orbits analogous to those of \cite{lus-spal}.

Let $P=MN$ be a parabolic subgroup of $H$ as above.  Let $X\in \sgo_M$ be a nilpotent element. Let  $\oc^M_X$ be the $M$-orbit of $X$. The  variety $\oc_X^M\oplus \sgo_N$ is irreductible. Since there are finitely many nilpotent orbits, there is a unique nilpotent $H$-orbit $\oc$ such that
  $$\oc\cap (\oc_X^M\oplus\sgo_N)$$
  is a Zariski  open dense subset of  $ \oc_X^M\oplus\sgo_N$. We denote $\oc$ by $I_P(X)$ and we shall call it \emph{the induced orbit}.

  \begin{proposition}\label{prop:Porb}
    The intersection  $I_P(X)\cap (\oc_X^M\oplus\sgo_N)$ is a single $P$-orbit.
  \end{proposition}

  \begin{preuve}     Let $Z\in \sgo_N$ such that   the element $Y=X+Z$ belongs to $I_P(X)\cap (\oc_X^M\oplus\sgo_N)$. Let  $\bc$ be  the variety of complete flags  of $E$-subspaces of $V$. As usual, we view $\bc$ as an $F$-variety. Let $\bc_Y$ be the subvariety of complete flags $V_\bullet$ such that $YV\subset V$.

It is not difficult to compute the dimensions of $\bc_Y$ and $H_Y$ in terms of the Jordan decomposition of $Y$ (one can use the method of \cite{Spal}). Then one can check the equality
\begin{align}
  \label{eq:dimHY}
\dim(H_Y)=2\dim(\bc_Y)+\dim(T).
\end{align}

    By using the group $B$ as a base point, we get an equivariant map $\pi_1:H \to \bc$. Let $\pi_2$ be the map from $H$ to the $H$-orbit $\oc_Y$ of $Y$ given by $h\mapsto h^{-1}Yh$. We have $\pi_2^{-1}(\oc_Y\cap \sgo_B)=\pi_1^{-1}(\bc_Y)$. The dimension of the fibers of $\pi_1$ and $\pi_2$ are respectively $\dim(B)$ and $\dim(H_Y)$. We deduce that
\begin{align}
  \label{eq:dimbase}
\dim(\oc_Y\cap \sgo_B)+\dim(H_Y)=\dim(\bc_Y)+\dim(B).
\end{align}
Combining \eqref{eq:dimHY} and \eqref{eq:dimbase}, we get
    \begin{align}
      \label{eq:dimNB}   \dim(\bc_Y)+\dim(\oc_Y\cap \sgo_B)=\dim(N_B).
    \end{align}
    We can also define a variety $\bc_X^M$ relative to $M$ and $X$. We have an identification $\bc_X^M\simeq \bc_Y^P=\bc_Y\cap\bc^P$ where $\bc^P$ is the variety of the complete flags that refine the flag associated to  $P$. In particular, $\dim(\bc_X^M)\leq \dim(\bc_Y)$. As before, one proves :
    $$\dim(\bc_X^M)+\dim(\oc_X^M\cap \sgo_{M\cap N_B})=\dim(N_B\cap M).$$
 Following the proof of theorem 1.3  in \cite{lus-spal}, we have
    \begin{align*}
      \dim(\bc_Y)+\dim(\oc_Y\cap \sgo_B)&\geq \dim(\bc_X^M)+\dim((\oc_X^M\oplus\sgo_N)\cap \sgo_B)\\
                                        &=\dim(\bc_X^M)+\dim(\oc_X^M\cap \sgo_{M\cap N_B})+\dim(\sgo_N)\\
                                        &=\dim(N_B\cap M)+\dim(\sgo_N)\\
      &=\dim(N_B).      
    \end{align*}
  By \eqref{eq:dimNB}, all the inequalities are in fact equalities. In particular, we deduce $\dim(\bc_Y)= \dim(\bc_X^M)$ and thus $\dim(H_Y)=\dim(M_X)$ where $M_X$ is the centralizer of $X$ in $M$.

    Let $\oc_Y^P$ be the $P$-orbit of $Y$ and $P_Y$ be its centralizer in $P$. We have 
    \begin{align*}
    \dim(I_P(X)\cap(\oc_X^M\oplus \sgo_N))\geq   \dim(\oc_Y^P)=\dim(P)-\dim(P_Y)&\geq \dim(P)-\dim(H_Y)\\
   &=\dim(P)-\dim(M_X)\\
      &=\dim(\oc_X^M\oplus \sgo_N)\\
&=  \dim(I_P(X)\cap(\oc_X^M\oplus \sgo_N)).
    \end{align*}
    Thus we have $\dim(\oc_Y^P)=\dim( I_P(X)\cap(  \oc_X^M\oplus \sgo_N))$. But one gets the same result for any $Y'\in  I_P(X)\cap (\oc_X^M\oplus\sgo_N)$. Thus the orbits $\oc_Y^P$ and $\oc_{Y'}^P$ must intersect by irreducibility of $I_P(X)\cap(  \oc_X^M\oplus \sgo_N)$.
  \end{preuve}
\end{paragr}

\begin{paragr} The following lemma is a variant of lemma 2.9.1 of \cite{cuft}.
  
  \begin{lemme}\label{lem:induite}
    There exists a finite family of polynomial maps $(\Phi_{i})_{i\in I}$  on $\sgo_{M}\oplus \sgo_{N}$ such that for any nilpotent  $M$-orbit $\oc$ in $\sgo_M$, there exists $I_\oc\subset I$ that satisfies the following  two properties:
\begin{enumerate}
\item For any $X\in \oc$ and  $Y\in   \sgo_N$, one has   $X+Y\in I_P(X)$ if and only if there exists $i\in I_\oc$ such that $P_{i}(X,Y)\not=0$. 
\item For any  $X\in \oc$, there exists $i\in I_\oc$ such that the restriction of $P_{i}(X,\cdot)$ to $\sgo_N$ is non-trivial.
\end{enumerate}
\end{lemme}

\begin{preuve}
Let $\oc$ be an nilpotent $M$-orbit in $\sgo_M$. Let $X\in \oc$ and $Y\in  \sgo_{N}$. Let $\oc'$ be the $P$-orbit of  $X$. By proposition \ref{prop:Porb}, we have $X+Y\in I_P(X)$ if and only if $\oc'$ is a dense open subset of $\oc\oplus \sgo_N$. The latter condition holds if and only if we  have the following equality between tangent spaces:
\begin{equation}
  \label{eq:CNS}
  [\pgo,X+Y]=[\mgo_P,X]\oplus\sgo_N.
\end{equation}
One always has the inclusion $[\pgo,X+Y]\subset[\mgo_P,X]\oplus\sgo_N$.  The dimension $d_\oc=\dim([\mgo_P,X]\oplus\sgo_N)$ does depend only on $\oc$. Thus the condition \eqref{eq:CNS} holds if and only if the rank of the $F$-map $\pgo \to \sgo_P$ given by $Z\mapsto [Z,X+Y]$ is at least $d_\oc$. This condition defines a dense Zariski  open subset  of $\sgo_P$ which does intersect non-trivially $X+\sgo_N$ since, for any nilpotent $X\in \sgo_M$, there exists $Y\in \sgo_N$ such that   $X+Y\in I_P(X)$. The statement is then clear.
\end{preuve}
\end{paragr}

\subsection{Nilpotent expansion and the homogeneity property}

\begin{paragr}
For any $F$-vector space $W$, we denote by $\Sc(W(\AAA))$ the space of complex Schwartz-Bruhat functions on $W(\AAA)$ and $\Cc(W(\AAA))\subset \Sc(W(\AAA))$  the subspace of smooth compactly supported functions. 
\end{paragr}

\begin{paragr}
Let $f\in \Sc(\sgo(\AAA))$. For any parabolic subgroup $P=MN$ of $H$ (see §\ref{S:parab}) and any $x\in H(\AAA)$, we define the constant term of $f$ by
$$f_{P,x}(X)=\int_{\sgo_N(\AAA)} f(\Ad(x)(X+U))\,dU
$$
where the Haar measure on $\sgo_N(\AAA)$ is normalized in such a way that the quotient  $\sgo_N(\AAA)/\sgo_N(F)$ equipped with the quotient measure by the counting measure is of volume $1$. This formula defines a function $f_{P,x}$ in the Schwartz-Bruhat space $\Sc(\sgo_M(\AAA))$.
\end{paragr}

\begin{paragr}\label{S:weak-cusp}
We say that $f$ is \emph{weakly cuspidal} if $f_{P,x}$ vanishes on the subset $\nc(\AAA)\cap \sgo_M(\AAA)$ for any proper  parabolic subgroup $P\subsetneq H$ and any $x\in H(\AAA)$. 
\end{paragr}

\begin{paragr}[A convergence result.] ---  Let $\nc(F)/H(F)$ be the finite set of $H(F)$-orbits on $\nc(F)$. Let $\oc\in \nc(F)/H(F)$ be a nilpotent orbit. For any  $f\in \Sc(\sgo(\AAA))$ and $h\in H(\AAA)$, we define
$$k_{\oc}(f,h)=\sum_{X\in \oc} f(\Ad(h^{-1})X).
$$

 \begin{theoreme}\label{cor:cv-inf}
    Let $f\in \Sc(\sgo(\AAA))$ be a  weakly cuspidal function. The integral
    \begin{align*}
     J_\oc(f)=  \int_{[H]^1}  k_{\oc}(f,h)\, dh
    \end{align*}
    is absolutely convergent.
  \end{theoreme}

The proof will be given in §\ref{S:proof-cv-inf} below. We have the following corollary which is also a simple consequence of a much more general result of Li (see theorem 1.1 of Li's paper \cite{Li2}). 

\begin{corollaire}\label{thm:cv}
  Let $f \in \Sc(\sgo(\AAA))$ be a weakly cuspidal function. The integral 
$$\int_{[H]^1} |\sum_{X\in \nc(F)} f(\Ad(h^{-1})X) |\,dh 
$$
is convergent.
\end{corollaire}

\end{paragr}

\begin{paragr}[Homogeneity of nilpotent integrals.] ---  Let $v$ be a place of $F$. For any  $f\in \Sc(\sgo(\AAA))$ and any $t\in F_v^\times$, let $f_t$ be the function in $\Sc(\AAA)$ defined by $f_t(X)=f(t^{-1}X)$.

 \begin{theoreme}\label{thm:homog2}
     Let $f\in \Cc(\sgo(\AAA))$ be a weakly cuspidal compactly supported function.  Let $v$ be a place of $F$. There exists a  bound $\eta >0$ such  for any  $t_0,t $ in $F_v^\times$ such that $|t_0|_v<\eta$ and  $|t|_v\leq 1$  we have
     \begin{align*}
       J_\oc(f_{tt_0})= |t|_v^{\dim(\oc)/2}J_{\oc}(f_{t_0}).
     \end{align*}
   \end{theoreme}

  \begin{remarque}
The proof is given in section \ref{ssec:refinedhomog} but it is built upon several preliminaries results. Indeed to get theorem \ref{thm:homog2}, the convergence of theorem\ref{cor:cv-inf}  is not sufficient. The main difficulty is that we cannot a priori intervert the integral and the sum over $\oc$: this would in general make appear an infinite volume related to a centralizer. We shall rather use a roundabout method inspired by \cite{suvta}. As a byproduct,  we shall obtain an expression of   $J_{\oc}(f_{t_0})$ in terms of an explicit zeta integral depending on $f$ (see theorem \ref{thm:comput}). Here we will need a mild assumption on the support of $f_{t_0}$ hence the bound on $t_0$.
  \end{remarque}
  
For future reference, let's state a simple corollary.

\begin{corollaire}\label{thm:homog}
Under the hypothesis of theorem \ref {thm:homog2}, for any  $t\in F_v^\times$ such that $|t|_v\leq 1$ we have
$$\int_{[H]^1} \sum_{X\in \nc(F)} f_{tt_0}(\Ad(h^{-1})X) \,dh =\sum_{\oc \in \nc(F)/H(F)} |t|_v^{\dim(\oc)/2} J_{\oc}(f_{t_0}).$$
In particular, when $t\to 0$ the expression above is equivalent to
$$J_{(0)}(f_{t_0})=\vol([H]^1) f(0)$$
where $(0)$ is the orbit of $0$.
  \end{corollaire}
  \end{paragr}

\subsection{Refined convergence results}\label{ssec:refinedcv}

\begin{paragr}\label{S:cv-inf}
  Let $\oc\in \nc(F)/H(F)$ be a nilpotent orbit. 
Let $P=MN$ be  standard parabolic subgroup  of $H$ (cf. §\ref{S:parab}). We borrow notations from \cite{suvta}; to do this, we identify $H$ with $GL(n,E)$ and $P$  with a standard parabolic subgroup of $GL(n,E)$ by the choice of the basis of §\ref{S:parab}.  We follow §§ 1.6 and 2.7 of  $\cite{suvta}$ (relatively to  the base field $E$). We have the Harish-Chandra map from $H_P$ from $H(\AAA)$ to some real vector space $a_P$. Restricted to $P(\AAA)$ this is a morphism and we denote by  $M(\AAA)^1 $ the intersection of $M(\AAA)$ with  the kernel of $H_P$. Let $\tau_P$ be the characteristic function of the acute Weyl chamber in $a_P$ associated to $P$.  Let $Z_P$ be the maximal central $F$-split torus in $M$ and let $A_P$ be the neutral component of the group of $\RR$-points of the split component of $\Res_{F/\QQ}(Z_P)$.  We have $M(\AAA)=M(\AAA)^1 A_P$. The function $F^P$ is the characteristic function of some compact subset of $A_PM(F)N(\AAA)\back H(\AAA)$. Note that both $F^P$ and $H_P$ depend on the choice of the compact subgroup $K$ of §\ref{S:maxcpt}. 

  \begin{theoreme}\label{thm:cv-inf}
Let $f\in \Sc(\sgo(\AAA))$ be a weakly cuspidal function.
    The integral 
$$\int_{P(F)\back H(\AAA)^1} F^P(h) \tau_P(H_P(h)) k_{\oc}(f,h)\, dh
$$
is absolutely convergent.
  \end{theoreme}

  The proof will be given in § \ref{S:preuve-cvinf} below. 

 \begin{remarque}
    Theorem \ref{thm:cv-inf} is  a simple analog of proposition 3.5.1 and corollary 3.2.2 of \cite{cuft}.
  \end{remarque}
\end{paragr}

\begin{paragr}[Proof of  theorem \ref{thm:cv-inf}.] ---\label{S:preuve-cvinf} It follows the lines of the proof of theorem 3.2.1 of  \cite{cuft}. For the reader's convenience, we will sketch the main steps and the simplications in our case. The case where $P=H$ is obvious since then $h\mapsto F^P( h) \tau_P(H_P( h))=F^H(h)$ is compactly supported on $[H]^1$. So we assume $P\subsetneq G$. We denote by $\Delta_P$ the set of simple roots of $Z_P$ in $N$ and $\rho_P$ the half-sum of roots of $Z_P$ on $N$.  It is easy to see that theorem \ref{thm:cv-inf} is a direct consequence of the following majorization.

\begin{lemme}\label{lem:1}
Let $f\in \Sc(\sgo(\AAA))$ be a weakly cuspidal function. Let $\Om$ be a compact subset of $N(\AAA)M(\AAA)^1K$.
There exist  $\eps>0$ and $c_0>0$ such that 
 $$ \exp(-\bg 2\rho_P,H_P(a)\bd)|k_{\oc}(f,ah)| \leq c_0 \cdot\prod_{\al\in \Delta_P}\al(a)^{-\eps}$$
for all $h\in \Om $ and all $a\in A_P$ such that $\tau_P(H_P(a))=1$.
\end{lemme}

Lemma \ref{lem:1} itself is a straightforward consequence of 

\begin{lemme}\label{lem:2}
  Let $f$ and $\Om$ be as in lemma \ref{lem:1}. Let $\al\in \Delta_P$. There exists $c_0>0$ such that 
 $$ \exp(-\bg 2\rho_P,H_P(a)\bd)|k_{\oc}(f,ah)| \leq c_0 \cdot \al(a)^{-1}$$
for all  $h\in \Om $ and all $a\in A_P$ such that $\tau_P(H(a))=1$.
\end{lemme}
\end{paragr}

\begin{paragr}[Proof of lemma \ref{lem:2}.] ---   The simple root  $\al\in \Delta_P$ defines a maximal parabolic subgroup $Q$ that contains $P$. We denote $Q=LR$ be the standard Levi decomposition of $Q$ where $R$ is the unipotent radical of $Q$. Let $\bar{R}$ be the unipotent radical of the opposite  parabolic subgroup.

Let's denote by $Y\mapsto Y_{\bar R}$ the projection of $\sgo$ on $\sgo_{\bar R}$ acording to the decomposition $\sgo=\sgo_R\oplus\sgo_L\oplus\sgo_{\bar R}$. To prove lemma \ref{lem:2}, we split the sum $\sum_{X\in \oc} f(\Ad(a h)^{-1}X)$ into the following three contributions:
\begin{align}
  \label{c:1}
\sum_{X\in \oc, X_{\bar R}\not=0} f(\Ad(a h)^{-1}X) \, ;
\end{align}

\begin{align}
  \label{c:2}
  \sum_{X\in \oc\cap \sgo_Q(F)} f(\Ad(a h)^{-1}X) - \sum_{X\in \sgo_L(F), I_Q(X)=\oc} \sum_{Y\in \sgo_{R}(F)}f(\Ad(a h)^{-1}(X+Y)) \,;
\end{align}

\begin{align}
  \label{c:3}
\sum_{X\in \sgo_L(F), I_Q(X)=\oc} \sum_{Y\in \sgo_R(F)}f(\Ad(a h)^{-1}(X+Y)).
\end{align}

For the contribution \eqref{c:1} we have a better majorization (see proof of lemma 3.8.2 in \cite{cuft}): for any integer $k\geq 1$, there exists $c_0>0$ such that 
 $$ \exp(-\bg 2\rho_P,H_P(a)\bd)|\sum_{X\in \oc, X_{\bar R}\not=0} f(\Ad(a h)^{-1}X)| \leq c_0 \cdot \al(a)^{-k}$$
for all $h\in \Om $ and all $a\in A_P$ such that $\tau_P(H(a))=1$. 
For the contribution \eqref{c:3}, we  have the same kind of majorization. But to see it, we need to introduce a non-trivial continuous additive character $F\back \AAA \to \CC^\times$ and a non-degenerate bilinear form on $\sgo$ given by $\bg X, Y\bd=\trace(XY)$ (this is the trace of an $F$-endomorphism of $V_E$). Then, using the Poisson summation formula for the sum over $\sgo_R(F)$ and the fact that $f$ is weakly cuspidal, one gets
\begin{align*}
  \sum_{X\in \sgo_L(F), I_Q(X)=\oc} \sum_{Y\in \sgo_R(F)}f(\Ad(a h)^{-1}(X+Y))\\
=\sum_{X\in \sgo_L(F), I_Q(X)=\oc} \sum_{Y\in \sgo_{\bar{R}}(F), Y\not=0} \int_{\sgo_R(\AAA)}   f(\Ad(a h)^{-1}(X+U))\psi(\bg Y, U\bd)\,dU.
\end{align*}
At this point we can conclude as in the proof of lemma 3.8.3 of \cite{cuft}. The most difficult contribution is thus \eqref{c:2} which can be written as
\begin{align*}
  \sum_{X\in \sgo_L(F), I_Q(X)\not=\oc}\sum_{Y\in \sgo_R(F), X+Y\in \oc} f(\Ad(a h)^{-1}X) \\
- \sum_{X\in \sgo_L(F), I_Q(X)=\oc} \sum_{Y\in \sgo_{R}(F), X+Y\notin \oc}f(\Ad(a h)^{-1}(X+Y)).
\end{align*}
But then the argument is that used in §3.12 of \cite{cuft} with lemma \ref{lem:induite} above playing the role of lemma 2.9.1 of \cite{cuft}.
 
\end{paragr}

\begin{paragr}[Proof of theorem \ref{cor:cv-inf}.] --- \label{S:proof-cv-inf} Theorem \ref{cor:cv-inf} is a straightforward consequence of theorem \ref{thm:cv-inf} and the equality for any $h\in H(\AAA)$
  \begin{align}\label{eq:HN}
    \sum_P \sum_{\delta \in P(F)\back H(F)}  F^P(\delta h) \tau_P(H_P(\delta h))=1
  \end{align}
  where the sum is over all standard parabolic subgroups  $P$ of $H$ (see proposition 2.5.1 of \cite{suvta}).  
  \end{paragr}

\subsection{A limit formula for nilpotent contributions}

\begin{paragr}
In this section, we will get an expression for $J_\oc(f)$ as the residue at $s=0$ of a function $J_\oc(s,f)$ when $f$ is a weakly cuspidal function . In the next section, under a mild condition on $f$, the function $J_\oc(s,f)$ is expressed in terms of a zeta integral. In our context, this is a simple analog of results in \cite{suvta}. The homogeneity of $J_\oc(f)$ is then an easy consequence of this result. Once again, we borrow notations from \cite{suvta}: we will denote by $E^H$ the function $E^{GL_E(n)}$ defined in \cite{suvta} § 3.2 eq. (3.2.1): it is a characteristic  function on $[H]$. 
\end{paragr}

\begin{paragr}
Let $f\in \Sc(\sgo(\AAA))$ and  $\oc\in\nc(F)/H(F)$ be a nilpotent orbit. 

\begin{theoreme}\label{thm:limit}
  Assume that $f$ is weakly cuspidal. Then the integral
  \begin{align*}
    J_\oc(s,f)=\int_{[H]} E^H(h) k_\oc(f,h) |\det(h)|^s\, dh
  \end{align*}
is absolutely convergent for $s\in \CC$ with $\Re(s)>0$. Moreover,
\begin{align*}
  \lim_{s\to 0^+}     s J_\oc(s,f)=J_\oc(f)
\end{align*}
 where $ \lim_{s\to 0^+}$ means that the limit is taken over complex numbers $s$ such that $\Re(s)>0$.
\end{theoreme}

\begin{preuve}
  It is a variation on the proof of theorem 4.6.1 of \cite{suvta}. In fact, the hypothesis ``weakly cuspidal'' makes the situation even simpler.

  One shows that  for a standard parabolic subgroup $P\subset H$ the integral
  $$  J_\oc^P(s,f)=\int_{P(F)\back H(\AAA)} E^H(h) F^P(h) \tau_P(H_P(h)) k_{\oc}(f,h) |\det(h)|^s\, dh
  $$
  is absolutely convergent for $s\in \CC$ with $\Re(s)>0$. The convergence is a straightforward consequence of theorem \ref{thm:cv-inf} and the fact that the function $E^H(h)$ has a simple expression when  $F^P(h) \tau_P(H_P(h))=1$, see expression (3.2.1) of \emph{ibid.}. Moreover (as in  the proof of proposition 4.5.1 of \emph{ibid.}), one has
  $$\lim_{s\to 0^+} s J_\oc^P(s,f)=\int_{P(F)\back H(\AAA)^1}F^P(h) \tau_P(H_P(h)) k_{\oc}(f,h) \, dh
  $$
 One  gets the theorem by adding the contributions of the various parabolic subgroups $P$, see \eqref{eq:HN}.
\end{preuve}
\end{paragr}

\subsection{Computation of a nilpotent integral}

\begin{paragr} \label{S:di}Let $X\in \nc(F)$. Using the proof of lemma \ref{lem:classif-orbit}, one can show that there exist
  \begin{itemize}
  \item a integer $r\geq 1$ 
  \item a decomposition
$$V_E= \bigoplus_{1\leq i \leq j\leq r} V_j^i$$
$d_j=\dim(V_j^i)$ does not depend on $i$;
\item a basis $(e^i_{k,j})_{1\leq k \leq d_j}$ of $V_j^i$;
\item in the basis   $(e^i_{k,j})_{1\leq i \leq j \leq r, 1\leq k \leq d_j}$ we have $$Xe_{k,j}^i=\left\lbrace
  \begin{array}{l}
    e_{k,j}^{i-1} \text{ si } i> 1 ;\\
0  \text{   si } i= 1.
  \end{array}\right.
$$
\end{itemize}
We may order the basis by  $e^i_{k,j}< e^{i'}_{k',j'}$ if and only if one the following conditions are satisfied.

\begin{itemize}
\item $i<i'$ ;
\item $i=i'$ et $j>j'$ ;
\item $i=i'$, $j=j'$ et $k<k'$.
\end{itemize}

In this basis, the matrix of $X$ is given by 
$$
\begin{pmatrix} 0_{d_1+\ldots +d_r} & \underset{0_{}}{I_{d_2+\ldots+d_r  }} & 0 & 0&  0 \\   &0_{d_2+\ldots+d_r}  & 0 & 0& 0\\ & & \ddots & \underset{0}{I_{d_{r-1}+d_r}} & 0\\ & & & 0_{d_{r-1}+d_r} & \underset{0}{I_{d_{r}}} \\ & & & & 0_{d_r}\end{pmatrix}.
$$

In the following we may replace $X$ by a conjugate. So we can and we shall assume that the ordered basis $(e^i_{k,j})$ is the basis chosen in §\ref{S:parab}. Thanks to this basis, we shall also identify the group $H(F)$ with $GL(n,E)$ and the $F$-space $\sgo$ with the space $\mathfrak{gl}(n,E)$ of $n\times n$ matrices with coefficients in $E$. Note that trough these identifications, the action $\Ad$ of the group $H(F)$  on $\sgo$ is  the $\sigma$-conjugation of $GL(n,E)$ on $\mathfrak{gl}(n,E)$ (given by $(y,Y)\mapsto y Y\sigma(y)^{-1}$).
\end{paragr}

\begin{paragr}
Let 
$$\pgo=\mgo\oplus \ngo$$
where
\begin{itemize}
\item $\mgo= \bigoplus_{1\leq i \leq j\leq r} \Hom_E(V_j^i,V_{j}^{i})$;
\item $\ngo=\bigoplus \Hom_E(V_j^i,V_{j'}^{i'})$ where the sum is over $1\leq i \leq j\leq r$ et $1\leq i' \leq j'\leq r$ such that $i>i'$ or $i=i'$ and $j<j'$.
\end{itemize}

Let $M$ and $N$ be the $F$-subgroups of $H$ of Lie algebras $\mgo$ and $\ngo$. We get a    parabolic   subgroup $P=MN$ with a Levi decomposition. Recall that  $H_X$ denotes the  centralizer of $X$ in $H$. Let $M_X=M\cap H_X$ and $N_X=N\cap H_X$. One has $H_X\subset P$ and $H_X=M_XN_X$ is a Levi decomposition. The groups $M_X$ and $M$ can be respectively identified to $\prod_{1\leq j\leq r} GL_E(V_j^j)$ and    $\prod_{1\leq i \leq j \leq r} GL_E(V_j^i)$. The inclusion $M_X\subset M$ is given by the diagonal embedding of $GL_E(V_j^j)$ in $\prod_{1\leq i \leq j } GL_E(V_j^i)$: by ``diagonal'', we mean that we identify  $GL_E(V_j^j)$ to  $GL_E(V_j^i)$  via the $F$-isomorphism $V_j^j\simeq V_j^i $ given by $X^{j-i}$.
Let
\begin{align*}
  \ugo_X=\sgo \cap \left( \bigoplus \Hom_F(V_j^i,V_{j'}^{i'})\right)
\end{align*}
 where the sum is over  $1< i \leq j\leq r$ and $1\leq i' \leq j'\leq r$ such that $i-1>i'$ or $i=i'+1$ and $j<j'$.

\begin{lemme}\label{lem:isom-orb}
  The map $n\mapsto nXn^{-1}-X$ induces an isomorphism from $N_X\back N$ onto $\ugo_X$.
\end{lemme}

\begin{preuve}
  It is similar to the proof of proposition 4.5.1 of \cite{scuft}.
\end{preuve}

\end{paragr}

\begin{paragr}[Zeta integral.] --- \label{S:zeta} Let $L= \prod_{1\leq i < j \leq r} GL_E(V^i_j)$. An element $A$ of  $L$ is written  $(A_{i,j})_{ 1\leq i <j\leq r}$ with $A_{i,j}\in GL_E(V_j^i)$. For any $A\in L(\AAA)$, let 
\begin{align*}
  \Delta_A=\begin{pmatrix} 0_{d_1+\ldots +d_r} & \underset{0_{}}{\Delta_1(A)} & 0 & 0&  0 \\   &0_{d_2+\ldots+d_r}  & 0 & 0& 0\\ & & \ddots & \underset{0}{\Delta_{r-2}(A)} & 0\\ & & & 0_{d_{r-1}+d_r} & \underset{0}{\Delta_{r-1}(A)} \\ & & & & 0_{d_r}\end{pmatrix}
\end{align*}
where  
\begin{align*}
   \Delta_i(A)= \begin{pmatrix} A_{i,r} &  &  & \\  & A_{i,r-1}  &  & \\ & & \ddots &  \\ & & & A_{i,i+1}
  \end{pmatrix}
\end{align*}
for any $1\leq i \leq r-1$.

For any $f\in \Sc(\sgo(\AAA))$ let's define
$$f_X^K (A)= \int_{\ugo_{X}(\AAA)} \int_K f(k^{-1}( U+\Delta_A)\sigma(k))\,dU dk$$
and for $s\in \CC$ such that $\Re(s)>0$  the zeta function
$$
Z_X(f,s)= \int_{L(\AAA)}  f^K_{X}(A)  \delta(A,s) dA 
$$
with
$$\delta(A,s)= \prod_{1\leq i < j \leq r}  |\det(A_{i,j})|^{d_i+\ldots+d_j+(j-i)s  }
$$
and $dA$ is a Haar measure on $L(\AAA)$ (normalized as in §\ref{S:Haar}). Recall that $\det$ is a shortcut for $N_{E/F}\circ\det$ (cf. §\ref{S:det}).

As in \cite{suvta} §8.3, the integral is convergent and defines a holomorphic   function on the domain $\Re(s)>0$.

\begin{lemme}\label{lem:homog}
Let $\oc$ be the $H$-orbit of $X$. Let $v$ be a place of $F$.   For any $t\in F_v^\times$, we have
\begin{align*}
  Z_X(f_t,s)=|t|_v^{\dim(\oc)/2+ c s}   Z_X(f,s)
\end{align*}
where $c=\sum_{1\leq j\leq r} j(j-1)d_j$.
\end{lemme}

\begin{remarque}
  As usual in this paper, $\dim(\oc)$ is the dimension of $\oc$ over $F$.
\end{remarque}

\begin{preuve}
Clearly we have the homogeneity property with the exponent:
\begin{align*}
    2 \sum_{1\leq i < j \leq r} d_j (d_i+\ldots+d_j+(j-i)s )+\dim(\ugo_X).
\end{align*}
We have on one hand
\begin{align*}
    2 \sum_{1\leq i < j \leq r} d_j (j-i)s=  \sum_{1< j \leq r} j(j-1)d_j s.
\end{align*}
On the other hand it suffices to show that
  \begin{align*}
    2 \sum_{1\leq i < j \leq r} d_j(d_i+\ldots+d_j)+\dim(\ugo_X)=\dim(\oc)/2.
  \end{align*}
One can compute $\dim(H_X)$ as follows
\begin{align*}
  \dim(H_X)/2&=\sum_{1\leq j,j'\leq r} d_j d_{j'}\min(j,j')\\
&=2 \sum_{1\leq i<j\leq r} i d_i d_{j}+ \sum_{1\leq j\leq r} jd_j^2\\
&= 2 \sum_{1\leq i < j \leq r} d_j(d_i+\ldots+d_j) -2 \sum_{1\leq j\leq r} (j-1)d_j^2+  \sum_{1\leq j\leq r} jd_j^2\\
&= 2 \sum_{1\leq i < j \leq r} d_j(d_i+\ldots+d_j) - \sum_{1\leq j\leq r} jd_j^2+2 \sum_{1\leq j\leq r} d_j^2\\
&= 2 \sum_{1\leq i < j \leq r} d_j(d_i+\ldots+d_j)  -\dim(M)/2+\dim(M_X).
\end{align*}
Using the equality $\dim(\ugo_X)=\dim(N/N_X)$ (cf. lemma \ref{lem:isom-orb}), we get
\begin{align*}
   2 \sum_{1\leq i < j \leq r} d_j(d_i+\ldots+d_j)+\dim(\ugo_X)&=  \dim(H_X)/2 + \dim(M)/2-\dim(M_X)+\dim(N)-\dim(N_X)\\
&= (\dim(M)+2\dim(N))/2 -\dim(H_X)/2\\
&=(\dim(H)-\dim(H_X))/2=\dim(\oc)/2.
\end{align*}
\end{preuve}
\end{paragr}

\begin{paragr}[Computation of the nilpotent integral.] ---  \label{S:computation}We will denote by $\theta_X(s)$ the function defined in \cite{suvta} § 7.2 relatively to the field $E$ and the datum $(d_1,\ldots,d_r)$ (cf. §\ref{S:di}). The main property of $\theta_X(s)$ we retain is that it  is holomorphic for $s\in \CC$ such that $\Re(s)>0$.

 \begin{theoreme}\label{thm:comput}
Let $f\in \Cc(\sgo(\AAA))$ and $v$ a place of $F$. There exists a  bound $\eta >0$ such that for any  $t$ in $F_v^\times$ such that $|t|_v<\eta$ we have 
\begin{align*}
  J_\oc(f_t,s)=\theta_X(s)\cdot Z_X(f_t,s).
\end{align*}
\end{theoreme}

  \begin{preuve} It is analogous to the proof of theorem 9.1.1 of \cite{suvta}. By the Iwasawa decomposition $H(\AAA)=P(\AAA) K$, we can write $g=mnk$ et $m=(m_j^i)\in M\simeq\prod_{1\leq i \leq j \leq r} GL_E(V_j^i)$, $n\in N(\AAA)$ and $k\in K$. Let $R=GL_E(d_1+\ldots+d_r)$ and  $r(g)\in R(\AAA)$ be the matrix extracted from the first  $d_1+\ldots+d_r$ rows and columns of $mn$.
As in \cite{suvta} p.115, one shows that $E^G(g)=1$ if and only if $E^{R}(r(g))=1$. The main new ingredient in our context is the observation that $E^G(g)=E^G(\sigma(g))$. Then the same kind of computations as those of \cite{suvta} p.116 leads to the statement (see also remark 9.1.2 of \emph{ibid.}).
\end{preuve}

\end{paragr}

\subsection{Proof of theorem \ref{thm:homog2}}\label{ssec:refinedhomog}

    Let $f\in \Cc(\sgo(\AAA)) $ a weakly cuspidal function. Recall that we have defined integrals $J_\oc(f,s)$ (see theorem \ref{thm:limit}). Let $\eta$ be the bound given by theorem \ref{thm:comput}. Let $t_0,t $ in $F_v^\times$ such that $|t_0|_v<\eta$ and  $|t|_v\leq 1$. By theorem \ref{thm:comput}, $J_\oc(f_{tt_0},s)$ and $J_\oc(f_{t_0},s).$ can be expressed in terms of a zeta function for which we have a homogeneity property (see lemma \ref{lem:homog}). We deduce that we have for the constant $c$ of lemma  \ref{lem:homog}
$$ J_\oc(f_{tt_0},s)= |t|_v^{\dim(\oc)/2+ c s}  J_\oc(f_{t_0},s).
$$
Taking the product with $s$ and then the limit on $s\to 0^+$ given by theorem  \ref{thm:limit}, we get the result.

\section{The unipotent contribution}\label{sec:unip}
  
\subsection{Algebraic situation}

\begin{paragr}
  We follow notations of section \ref{ssec:algprel}. Let $D$ be a quaternion algebra over $F$ equipped with a fixed embedding $E\hookrightarrow D$. Note that $D$ may be split. Then $V_D=V\otimes_E D$ is a  right $D$-module. Let $G=\Aut_D(V_D)$ viewed as an $F$-group. Let $\eps\in G$ given by the left multiplication by $\sqrt{\tau}$. Let $\theta$ the involution of $\End_F(V_D)$ given by $\theta(X)=\eps X \eps^{-1}$.   Let $H'\subset G$ be the subgroup fixed by $\theta$.  Let $S'\subset G$ be the $F$-variety of automorphism $g\in  \Aut_D(V_D)$ such that $g\eps=\eps g^{-1}$. The map
  \begin{align}\label{eq:rho}
    \rho:g\mapsto g\theta(g)^{-1}
  \end{align}
  induces an isomorphism from $G/H'$ onto $S'$ (cf. \cite{GuoPJ}). The action of $G$ by left translations on $G/H'$ gives an action of $G$ by $\theta$-conjugation for which $\rho$ is equivariant: we have $\rho(g_1 g_2)= g_1 \rho(g_2)\theta(g_1)^{-1}$.  This action induces an action by conjugation of $H'$ on $S'$. 
\end{paragr}

\begin{paragr}
  Let $\uc\subset G$ be the unipotent variety. We define $\uc_{S'}=\uc\cap S'$ and $\uc_G=\rho^{-1}(\uc_{S'})$. One knows (\cite{GuoPJ} lemma 3.2) that $\uc_G(F)=H'(F)\uc_{S'}(F) H'(F)$.
\end{paragr}

\begin{paragr}
  The tangent space of $S'$ at $\Id_{V_D}$  is denoted by $\sgo'$ : it is the space of $X\in \End_D(V_D) $ such that $X\eps+\eps X=0$. Let $\nc_V$ be the cone of nilpotent elements in $\End_D(V_D)$. Let $\nc'=\nc_V\cap \sgo'$. The usual exponential map denoted by $\exp$ induces an isomorphism from $\nc_V$ to $\uc$ and also from $\nc'$ to $\uc_{S'}$.

\end{paragr}

\begin{paragr}
  The map $\End_E(V)\to \End_D(V_D)$ given by $\varphi\mapsto \varphi\otimes \Id_D$ gives an identification of $H$ with $H'$ and $\sgo$ with $\sgo'$ and $\nc$ with $\nc'$. Hence we can freely use all the results of section \ref{sec:inf} for the action of $H'$ on $\sgo'$.

  To simplify the notations, we will suppress the superscript $'$ and we will not distinguish between $H$ and $H'$, $S$ and $S'$ etc. For $v$ a place of $F$, the measures used on the groups of $F_v$-points are Haar measures (we do not need any normalization).
\end{paragr}

\subsection{Main results}

\begin{paragr}
  Let $v$ be a place of $F$. For any  $\Phi\in \Cc(G(F_v))$, there is a unique function, denoted by $\Phi_{S}\in  \Cc(S(F_v))$ such that 
$$\Phi_{S}( \rho(x))=\int_{H(F_v)} \Phi(xh)\,dh.$$
The map $\Phi\mapsto \Phi_S$ is a surjection from $\Cc(G(F_v))$ onto  $\Cc(S(F_v))$.
\end{paragr}

\begin{paragr}[Very cuspidal test functions.] --- We shall say that $\Phi\in \Cc(G(F_v))$ is \emph{very cuspidal} if one has 
$$\int_{N(F_v)}\Phi(x ny)\, dn=0$$
for any parabolic subgroup $P\subsetneq G$ and any  $x,y\in G(F_v)$. Here $N$ is the unipotent radical of $P$. 

\begin{remarque}\label{rq:mat-coef}
Assume $v$ is finite.  Let $\tilde{\Phi}$ be  a matrix coefficient of a supercuspidal representation of $G(F_v)/Z(F_v)$ where $Z$ is the center of $G$. Let $\Phi\in \Cc(G(F_v))$  such that 
  \begin{align*}
    \int_{Z(F_v)}\Phi(gz)\,dz=\tilde{\Phi}(g)
  \end{align*}
for any $g\in G(F_v)/Z(F_v)$ . Then $\Phi$ is very cuspidal.
\end{remarque}

Note that any very cuspidal function $\Phi\in \Cc(G(F_v))$ is such that
\begin{align}\label{eq:vanish}
\int_{N(F_v)/N_{H}(F_v)}\Phi_{S}(\rho(xn))\, dn=0
\end{align}
where $N_{H}=N\cap {H}$ and $dn$ is the quotient of the Haar measures on $N(F_v)$ and $N_{H}(F_v)$.

\end{paragr}

\begin{paragr}[Global setting.] --- We also have  a surjective map $\Phi\mapsto \Phi_{S}$ from $\Cc(G(\AAA))$ onto $\Cc(S(\AAA))$ given by
$$\Phi_{S}( \rho(x))=\int_{H(\AAA)} \Phi(xh)\,dh$$
for any $x\in G(\AAA)$.

  We shall say that  $\Phi\in \Cc(G(\AAA))$ is\emph{ very cuspidal} if there is a place $v$ such that if one writes $\AAA=F_v\times\AAA^v$, one has $\Phi=\Phi_v \otimes \Phi^v$ with $\Phi^v\in \Cc(G(\AAA^v))$ and $\Phi_v  \in \Cc(G(F_v))$ is  very cuspidal.

\end{paragr}

\begin{paragr}
  Let's define for $\Phi\in \Cc(G(\AAA))$ and $x,y\in G(\AAA)$
$$K_{\uc_G}(\Phi,x,y)=\sum_{\gamma\in \uc_G(F)} \Phi(x^{-1}\gamma y)$$
and for $h\in H(\AAA)$
$$K_{\uc_S}(\Phi,h)=\sum_{\gamma\in \uc_S(F)} \Phi_S(h^{-1}\gamma h ).$$
We have the simple relation for $h\in H(\AAA)$
\begin{equation}
  \label{eq:sym-ker}
  \int_{[H]}  K_{\uc_G}(\Phi,x,h)\,dx= K_{\uc_S}(\Phi,h).
\end{equation}

\begin{theoreme}\label{thm:homog-gp}
Let   $\Phi\in \Cc(G(\AAA))$ be  a very  cuspidal function. Let  $\oc\in \nc(F)/H(F)$ and let $v$ be a place of $F$. 
\begin{enumerate}
\item For any $t\in F_v^\times$,  the integral 
$$J_\oc^t(\Phi)=\int_{[H]^1} \sum_{X\in \oc} \Phi_S(h^{-1}\exp(t^{-1}X)h )\,dh$$
is absolutely convergent. For $t=1$, $J_\oc^1(\Phi)$ is denoted by $J_\oc(\Phi)$.
\item (Fine expansion) We have
$$\int_{[H]^1} K_{\uc_S}(\Phi,h)\, dh=\sum_{\oc\in \nc(F)/H(F)} J_\oc(\Phi)$$
where the left-hand side is absolutely convergent.
\item \label{ass:2}There exists a  bound $\eta >0$ such that for any $t$ in $F_v^\times$ such that $|t|_v<\eta$, we have
$$ J_\oc^t(\Phi)=\lim_{s\to 0^+} s \theta_X(s)  Z_X(f_t,s)$$
where $X\in \oc$ is the element considered in §\ref{S:di}, $\theta_X(s)$ (defined in §\ref{S:computation}) is holomorphic for $\Re(s)>0$ and does not depend on $\Phi$, the zeta function $Z_X(f_t,s)$ is defined in §\ref{S:zeta} relatively to any function $f_t\in \Cc(\sgo(\AAA))$ such that $f_t(Y)=\Phi_S(\exp(t^{-1}Y))$.
\item Let $\eta$ be the bound of \ref{ass:2}. For any  $t_0,t $ in $F_v^\times$ such that $|t_0|_v<\eta$ and  $|t|_v\leq 1$  we have
     \begin{align*}
       J_\oc^{tt_0}(\Phi)= |t|_v^{\dim(\oc)/2}J_{\oc}^{t_0}(\Phi).
     \end{align*}
   \end{enumerate}
 \end{theoreme}
 
\begin{preuve}
   Theorem \ref{thm:homog-gp} is deduced from similar results on $\sgo$ and the (standard) descent procedure to $\sgo$ that is explained in §\ref{S:proof-ass3}. More precisely, assertion 1 results from theorem \ref{cor:cv-inf}. Assertion 2 is a consequence of assertion 1 and the finiteness of nilpotent orbits. Assertion 3 is a combination of theorem \ref{thm:limit} and \ref{thm:comput}. Finally assertion 4 is a consequence of assertion 3 (see section \ref{ssec:refinedhomog}). 
 \end{preuve}
 
Let's state a corollary which is a straightforward consequence of theorem \ref{thm:homog-gp}.

\begin{corollaire} \label{cor:homog-gp}We use notations of theorem \ref{thm:homog-gp}. 
For any  $t_0,t $ in $F_v^\times$ such that $|t_0|_v<\eta$ and  $|t|_v\leq 1$  we have
\begin{align*}
  \int_{[H]^1} |\sum_{X\in \nc(F)} \Phi_S(h^{-1}\exp((tt_0)^{-1}X)h )|\, dh<\infty
\end{align*}
and 
$$\int_{[H]^1} \sum_{X\in \nc(F)} \Phi_S(h^{-1}\exp((tt_0)^{-1}X)h )\, dh=\sum_{\oc \in \nc(F)/H(F)} |t|_v^{\dim(\oc)/2} J_\oc^{t_0}(\Phi).$$
In particular, when $t\in F_v^\times $ goes to $0$, the expression 
$$\int_{[H]^1} \sum_{X\in \nc(F)} \Phi_S(h^{-1}\exp(t^{-1}X)h )\, dh
$$
is equivalent to 
$$\vol([H]^1) \int_{H(\AAA)}\Phi(h)\,dh.$$
\end{corollaire}

\begin{remarque}
For  $\Phi$ is very cuspidal, we have the equality
  $$
\int_{[H]^1} \int_{[H]} K_{\uc_G}(\Phi,x,y)\, dx dy=\int_{[H]^1} K_{\uc_S}(\Phi,h)\, dh.
$$
where the left-hand side is at least conditionnally  convergent.  One can prove in fact that it is absolutely convergent using mixed truncation operators (in the sense of the seminal paper \cite{JLR}) built upon the combinatorics of \cite{Li2}. 
\end{remarque}
\end{paragr}

\begin{paragr}[Descent to the tangent space.] --- \label{S:proof-ass3} Let $\vc_0$ be a finite set of places of $F$ containing the archimedean places and  a fixed place denoted by $v_0$. Let $A\subset F$ the ring of integers outside $\vc_0$. We assume that $\vc_0$ is large enough such that  all objects  $G$, $H$, $S$, $\uc_{S}, \nc$   come naturally from $A$-schemes by base change. We assume also that the exponential (denoted by $\exp$)  induces an isomorphism of $A$-scheme from $\nc$ to $\uc_{S}$. For $v\notin \vc_0$, let $\oc_v\subset F_v$ be the ring of integers.

  Let
$$\Phi=\Phi_{0}\otimes \Phi_{1} \otimes \Phi_2$$
where  $\Phi_{0}\in \Cc(G(F_0))$ (with $F_0=F_{v_0}$) is a very cuspidal function,  $\Phi_{1}$ is a test function on 
$$F_1=\prod_{ v\in \vc_0\setminus\{v_0\} } F_v$$ 
and $\Phi_2$ is the characteristic function of $\prod_{v\notin \vc_0} G(\oc_v)$.

Let $c:\sgo\to \cgo=\sgo//H$ be the categorical quotient. For all $v\in \vc_0$, we fix an open subset $\om^\flat_v\subset \cgo(F_v)$ containing $c(0)$ such that the exponential map $\exp$ is well-defined and induces an analytic diffeomorphism from $\om_v=c^{-1}(\om^\flat_v)$ onto an  open subset $\Om_v\subset S(F_v)$. Let $\zeta_v$ be a smooth function on $\cgo(F_v)$ with compact support included in $\om_v^\flat$ and with value $1$ in a neighborhood of $c(0)$.

Let's define $\zeta_0=\zeta_{v_0}$ and  $\om_0=\om_{v_0}$ and also $\zeta_1=\prod_{ v\in \vc_0\setminus\{v_0\} } \zeta_v$ and $\om_1=\prod_{ v\in \vc_0\setminus\{v_0\} } \om_v$.  We define functions $f_0\in \Cc(\sgo(F_{0}))$ and $f_1\in \Cc(\sgo(F_1))$ by: 
  \begin{align*}
    f_i(X)=\left\lbrace
    \begin{array}{l} \zeta_i(X) \Phi_{i, S}(\exp(X)) \text{  if  } X\in \om_i;\\
      0 \text{ otherwise};
    \end{array}\right.
  \end{align*}
  for any $i\in \{0,1\}$ and $X\in \sgo(F_i)$.
Let $f=f_0\otimes f_1 \otimes f_2\in \Cc(\sgo(\AAA))$ where $f_2$ is the characteristic function of $\prod_{v\notin \vc_0}\sgo(\oc_v)$.

We have for all $X\in \nc(\AAA)$ 
\begin{align*}
  \Phi_S(h^{-1}\exp(X)h )= f(\Ad(h^{-1})(X)).
\end{align*}

To complete the proof of theorem \ref{thm:homog-gp} it suffices to check the next lemma.

\begin{lemme}\label{lem:wcusp}
  The function $f$ is weakly cuspidal (in the sense of §\ref{S:weak-cusp}).
\end{lemme}

\begin{preuve}
   Let $X\in \sgo_M(F_0)\cap \nc(F_0)$. Clearly, it suffices to prove that
$$\int_{\sgo_N(F_0)} f_0(\Ad(x)(X+U))\,dU=0$$
for any $x\in H(F_v)$. But  for any $U\in \sgo_N(F_0)$ we have
\begin{align*}
   f_0(\Ad(x)(X+U))&=  \Phi_{0, S}(x\exp(X+U)x^{-1}).
\end{align*}
Let's write
\begin{align*}
    \exp(X+U) = \exp(X/2)   ( \exp(-X/2)\exp(X+U)\exp(-X/2)  ) \theta(\exp(X/2))^{-1}.
\end{align*}
We observe that $U\mapsto \exp(-X/2)\exp(X+U)\exp(-X/2)$ induces an isomorphism from $\sgo_N$ onto $S\cap N$. But the map $\rho$ of \eqref{eq:rho} induces an isomorphism from $N/N_{H}$ onto $S\cap N$. We get a bijection from $N(F_0)/N_H(F_0)$ onto $\sgo_N(F_0)$. By a change of variables, we get (up to a constant $c\not=0$) 
\begin{align*}
  \int_{\sgo_N(F_0)} f_0(\Ad(x)(X+U))\,dU=c \int_{N(F_0)/N_H(F_0)}  \Phi_{0, S}( \rho(x\exp(X/2) n))\, dn=0
\end{align*}
by the vanishing condition \eqref{eq:vanish}.
\end{preuve}
\end{paragr}

\bibliography{biblio}
\bibliographystyle{alpha}

\begin{flushleft}
Pierre-Henri Chaudouard \\
Université Paris Diderot (Paris 7)\\
 Institut de Mathématiques de Jussieu-Paris Rive Gauche \\
 UMR 7586 \\
 Bâtiment Sophie Germain \\
 Case 7012 \\
 F-75205 PARIS Cedex 13 \\
 France
\medskip

email:\\
Pierre-Henri.Chaudouard@imj-prg.fr \\

\end{flushleft}
\end{document}